\documentclass[preprint,12pt,leqno]{elsarticle}
\pdfoutput=1
\usepackage{preamble}
\journal{the Journal of Algebra}
\date{15th April 2025}
\begin{document}
\begin{frontmatter}
\title{A 2-dimensional torsion theory on symmetric monoidal categories}
\author[inst1]{Mariano Messora}
\ead{mariano.messora@unimi.it}
\affiliation[inst1]{organization={Dipartimento di Matematica, Università degli Studi di Milano},
            addressline={Via C. Saldini 50}, 
            city={Milano},
            postcode={20133}, 
            country={Italy}}
\begin{abstract}
In this paper, we describe a homotopy torsion theory on the category of small symmetric monoidal categories. By using natural isomorphisms as the basis for the nullhomotopy structure, this homotopy torsion theory exhibits some interesting 2-dimensional properties which could be the foundation for a definition of ``2-dimensional torsion theory''.

We choose symmetric 2-groups as torsion objects, thereby generalising a known pointed torsion theory in the category of commutative monoids where abelian groups are taken as torsion objects. In the final part of the paper we  carry out an analogous generalisation for the classical torsion theory in the category of abelian groups given by torsion and torsion-free abelian groups.
\end{abstract}
\begin{extrainfo}
    \copyright{} 2025. This manuscript version is made available under the CC-BY-NC-ND 4.0 license \url{https://creativecommons.org/licenses/by-nc-nd/4.0/}
\end{extrainfo}
\end{frontmatter}
\setcounter{tocdepth}{1}
\renewcommand{\baselinestretch}{0.95}\normalsize
\tableofcontents
\renewcommand{\baselinestretch}{1.0}\normalsize

\section*{Introduction}
\addcontentsline{toc}{section}{Introduction}
\label{Intro}
Torsion theories were originally introduced in the context of abelian categories (\cite{DICKSON}) as a generalisation of a well-known property of the category of abelian groups: this category includes two special classes of objects---torsion and torsion-free groups---which exhibit the following key features. 
\begin{enumerate}
\item \label{ab1} Any group homomorphism from a torsion group to a torsion-free group is necessarily zero;
\item \label{ab2} any abelian group $A$ fits into a short exact sequence of abelian groups $T\to A\to F$, where $T$ is a torsion group and $F$ is a torsion-free group.
\end{enumerate}
More recently, many different authors have studied wider generalisations of this notion, both in pointed and in non-pointed categories (see, for example, \cite{BOURN,CLEMENTINO,ROSICKY,JANELIDZE,GRANDIS2,FACCHINI}).
The objective of the present work is to introduce and study a torsion theory in a 2-dimensional context. One can already find a very interesting case of a torsion theory developed in the 2-category $\catCat$ of small categories in the recent work of Borceux, Campanini, Gran and Tholen (\cite{BORCEUX}); however, in that work, due to the use of the formalism of pretorsion theories, the 2-category $\catCat$ is essentially treated as a 1-category, and the 2-dimensional aspects are not fully taken into account. Our goal, instead, is to construct a torsion theory in a 2-category in such a way that the 2-dimensional structure is extensively incorporated into the theory. We achieve this goal by applying the formalism of \emph{homotopy torsion theories}, recently introduced in \cite{HTT}.

Our starting point is a pointed torsion theory on the (1-)category $\CMon$ of commutative monoids, first introduced by Facchini in \cite[Section~1.3]{FACCHINI2}. 
Torsion objects are given by abelian groups, seen as commutative monoids in which every element is invertible, while torsion-free objects are those monoids in which the only invertible element is the unit. In the literature, monoids with this property have been referred to as \emph{pure}, \emph{sharp} or \emph{reduced}; here, we adopt the first name. It is shown in \cite{FACCHINI2} that abelian groups and pure commutative monoids form indeed a torsion theory on $\CMon$.

The main objective of this paper is to generalise this torsion theory, extending it from the 1-dimensional context of commutative monoids to the 2-dimensional one of symmetric monoidal categories. 
However, this generalisation is not straightforward.
First, the 1-category of symmetric monoidal categories (with symmetric monoidal functors) is \emph{not} a pointed category, and therefore zero morphisms and short exact sequences cannot be defined in the traditional sense. 
Second, we know that symmetric monoidal categories actually form (with symmetric monoidal functors and monoidal natural transformations) a \emph{2-category}, and, as mentioned above, we wish to take this richer structure into account. 

The framework of homotopy torsion theories allows us to address both the above issues. In this setting, we prove that symmetric 2-groups are the torsion objects of a homotopy torsion theory in the category of symmetric monoidal categories, using particular natural isomorphisms as nullhomotopies. We also show that, thanks to the use of 2-cells as nullhomotopies, this homotopy torsion theory exhibits some remarkable properties in dimension 2 which may be the basis for a notion of ``2-dimensional torsion theory'' in future work. We conclude the paper by carrying out an analogous generalisation for the classical torsion theory on abelian groups given by torsion and torsion-free groups, thus describing another homotopy torsion theory which enjoys the same 2-dimensional properties.
\section{Preliminaries}
\label{Pre}
In this section, we recall some fundamental notions and results that will be needed for the construction of the new torsion theory. 
\subsection*{Monoidal categories and monoidal functors}
\label{MonCat}
Since the core idea is to transition from the 1-dimensional context of monoids to the 2-dimensional context of monoidal categories, we first establish the notation and conventions that will be used for monoidal categories.
For the basic theory of monoidal categories we refer for example to \cite{MACLANE} or \cite{BORCEUX2}; see also \cite{KELLY} and \cite{JOYAL} for coherence issues.
We denote the full structure of a monoidal category by
\begin{equation*}
	(\cat M, \ten_{\cat M}, I_{\cat M}, a_{\cat M}, \ell_{\cat M}, r_{\cat M}), 
\end{equation*}
where $\cat M$ is the underlying category, $\ten_{\cat M}$ is the tensor, $I_{\cat M}$ is the tensor unit, $a_{\cat M}$, $\ell_{\cat M}$, $r_{\cat M}$ are the associator, the left and the right unitors respectively. We will often simply write $(\cat M, \ten, I, a, \ell, r)$, or even just $\cat M$. When $\cat M$ is symmetric, we will denote the symmetric braiding by $b_{\cat M}$ (or simply $b$).

For a monoidal functor
\begin{equation*}
	F\colon \cat M\to\cat N,
\end{equation*}
we denote by $e_F$ and $m_F$ (or simply $e$ and $m$) the maps encoding coherence with the unit and the tensor product, respectively.
We always assume these maps to be isomorphisms.

We denote by $\MC$ the 2-category of small monoidal categories, monoidal functors and monoidal natural transformations,
and by $\SMC$ the 2-category of small symmetric monoidal categories, symmetric monoidal functors and monoidal natural transformations. 
As with any strict 2-category, we will sometimes implicitly regard $\MC$ and $\SMC$ as ordinary 1-categories, with no risk of confusion.

In what follows, the symmetry  $b$, the structure maps $a$, $\ell$ and $r$ of monoidal categories, as well as the coherence maps $e$ and $m$ of monoidal functors will often be omitted.
\subsection*{2-groups}
\label{2G}
As we generalise the torsion theory considered in \intro{}, replacing commutative monoids with symmetric monoidal categories naturally suggests replacing abelian groups with symmetric 2-groups. 

While the notion of (weak) 2-group has appeared in the literature before, it is in Ho\`ang Xu\^an S\'inh's doctoral thesis (\cite{SINH}) that we find a first systematic study of the concept, under the name of ``\emph{gr-cat\'egories}'' (see also \cite{BAEZ-SINH}). For convenience, we now recall the definition of 2-group (weak 2-group in \cite{BAEZ}).
\begin{definition}
\label{def-2G}
In a monoidal category, we say that an object $A$ is \emph{weakly invertible} if there exists an object $B$ (called a \emph{weak inverse of $A$}) such that
\begin{equation*}
	A\ten B\iso B\ten A\iso I.
\end{equation*}
A \emph{2-group} is a monoidal category $\cat G$ such that
\begin{itemize}
	\item every morphism in $\cat G$ is invertible (i.e.\ $\cat G$ is a groupoid);
	\item every object in $\cat G$ is weakly invertible.
\end{itemize}
A \emph{symmetric 2-group} is a 2-group which is symmetric as a monoidal category.
We denote by $\TG$ (respectively, $\catGr$) the 2-category of small 2-groups (respectively, small symmetric 2-groups), monoidal functors (respectively, symmetric monoidal functors) and monoidal natural transformations.
\end{definition}
The following elementary property that will be useful later (cf.\ for example \cite{BAEZ}).
\begin{proposition}
\label{basic-2g}
In a 2-group $\cat G$, for any object $A$, the endofunctors $A\ten -$ and $-\ten A$ of $\cat G$ are equivalences (with inverses $B\ten-$ and $-\ten B$ respectively, for any weak inverse $B$ of $A$). Furthermore, for each object $A$, it is always possible to find an object $B$ and isomorphisms
\begin{equation*}
	h\colon I\to A\ten B, \quad\quad k\colon B\ten A\to I,
\end{equation*}
such that $(A\ten-,B\ten -,h,k)$ forms an adjoint equivalence, i.e., it makes the following compositions identities:
\begin{equation*}
    \begin{tikzcd}[column sep =1.3em,row sep=1.5em]
        A\arrow[r]
        &I\ten A\arrow[r,"h\ten A"{, rotate=15, shift={(.2em,.1em)}}]
        &[.5em](A\ten B)\ten A\arrow[r]
        &A\ten(B\ten A)\arrow[r,"A\ten k"{, rotate=15, shift={(.2em,.1em)}}]
        &[.5em]A\ten I\arrow[r]
        &A,
        \\
        B\arrow[r]
        &B\ten I\arrow[r,"B\ten h"{, rotate=15, shift={(.2em,.1em)}}]
        &[.5em]B\ten(A\ten B)\arrow[r]
        &(B\ten A)\ten B \arrow[r,"k\ten B"{, rotate=15, shift={(.2em,.1em)}}]
        &[.5em]I\ten B\arrow[r]
        &B
    \end{tikzcd}
\end{equation*}
(the unnamed arrows are components of the unitors, of the associator or of their inverses).
Clearly $(B\ten-, A\ten-, k^{-1}, h^{-1})$ is an adjoint equivalence as well.
\end{proposition}

There is a natural way to assign a 2-group to any monoidal category: by considering its weakly invertible objects and isomorphisms between them. 
We recall this fact in the following definition.
\begin{definition}
\label{Pic-def}
Given any monoidal category $\cat M$, the \emph{Picard 2-group} of $\cat M$, which we denote by $\PIC M$, is the subcategory of $\cat M$ consisting of all the weakly invertible objects of $\cat M$ and isomorphisms between them. This inherits the monoidal structure of $\cat M$, and is therefore a 2-group. 
\end{definition}
(The Picard \emph{2-group} should not be confused with the Picard \emph{group}, which denotes the group of isomorphism classes of invertible objects and is briefly mentioned in \cref{HomInv}.)

To conclude the section, we state the following well-known property for future reference.
\begin{remark}
\label{2g-coref} In this remark, we regard $\TG$, $\catGr$, $\MC$ and $\SMC$ as 1-categories, as we are currently only interested in their 1-dimensional properties. 
$\TG$ is a coreflective subcategory of $\MC$, and $\catGr$ is a coreflective subcategory of $\SMC$. In both cases the coreflector is given by extending to 1-cells the mapping $\cat M\mapsto \PIC M$.
\end{remark}
In the next section we define a particular class of  monoidal categories whose properties are, in a sense, diametrically opposed to those of 2-groups.
\section{Purely monoidal categories}
\label{PMC}
In this section, we introduce and study the notion of \emph{purely monoidal category}, which is meant to provide a 2-dimensional counterpart to the notion of pure monoid (see \intro{}).
Consider the following equivalent statements about monoidal categories. 
\begin{proposition}
\label{Pure-eqdef}
Let $\cat M$ be a monoidal category, and recall that $\GP M$ is its Picard 2-group (see \cref{Pic-def}). The following are equivalent.
\begin{enumerate}[(1)]
	\item \label{Pure-eq}$\PIC M$ is monoidally equivalent to $\catone$, the monoidal category with a single object $I_\catone$ and a single morphism. 
	\item \label{Pure-id}All weakly invertible objects of $\cat M$ are isomorphic to $\I M$, and the only automorphism of $\I M$ is the identity. 
	\item \label{Pure-un}For any two weakly invertible objects in $\cat M$, there exists exactly one isomorphism between them.
\end{enumerate}
(We recall that a monoidal equivalence is simply an equivalence in the 2-category $\MC$.)
\end{proposition}
\begin{proof}
    The equivalence of \eqref{Pure-id} and \eqref{Pure-un} is obvious by the fact that $\I M$ is weakly invertible itself.

    Now call $\cat G$ the Picard 2-group of $\cat M$, and call $F$ the only (monoidal) functor from $\cat G$ to $\catone$. If \eqref{Pure-eq} holds, $F$ has an inverse $G\colon \catone\to \cat G$ with a monoidal natural isomorphism $\eta\colon GF\tto\id{\cat G}$. For any $A$ in $\cat G$, the composition $\comp{\eta_A^{-1}}{e_G}$ gives an isomorphism between $A$ and $\I M$ (where $e_G$ is the coherence map of $G$ described in \cref{MonCat}). By the naturality of $\eta$, this is the only map in $\cat G$ from $A$ to $\I M$. 
    Vice versa, if \eqref{Pure-un} holds, call $G\colon \catone\to\cat G$ the functor defined by $G I_\catone=\I M$. Plainly, this is a monoidal functor with the obvious coherence maps, and clearly $FG=\id\catone$. For every object $A$ in $\cat G$, the unique isomorphism between $A$ and $\I M$ provides an isomorphism between $GFA$ and $A$. The fact that these isomorphisms are the components of a monoidal natural isomorphism from $GF$ to $\id{\cat M}$ follows from the uniqueness hypothesis in \eqref{Pure-un}.
\end{proof}
\begin{remark}
Observe that \cref{Pure-eqdef} can be specialised to the category $\SMC$ simply by noticing that, if $\cat M$ is a symmetric monoidal category and $\PIC M$ is equivalent to $\catone$ in $\MC$, then $\PIC M$ is equivalent to $\catone$ in $\SMC$. 
\end{remark}
\begin{definition}
\label{def-Pu}
We say that a monoidal category $\cat M$ is \emph{pure} (or that $\cat M$ is a \emph{purely monoidal category}) if any (and hence all) of the conditions in  the above proposition hold. When a purely monoidal category is symmetric as a monoidal category, then we say that it is a \emph{symmetric purely monoidal category}. We denote by $\catPMC$ (respectively, $\catPu$) the 2-category of small purely monoidal categories (respectively, small symmetric purely monoidal categories), monoidal functors (respectively, symmetric monoidal functors) and monoidal natural transformations.
\end{definition}
\begin{remark}
    The requirement of a \emph{unique} isomorphism in the definition of purely monoidal category can be intuitively understood as a 2-dimensional analogue of the requirement of a unique invertible element in the definition of pure monoid (see the last part of \cref{prop-hom-inv}).

    Moreover, the choice of 2-groups as torsion objects (and of the nullhomotopy structure from \cref{Null-struct}) uniquely determines the torsion-free objects (see \cref{CNH} and in \cite[Proposition~7.3]{HTT}), thereby forcing this definition of purely monoidal category.
\end{remark}
\begin{example}
Cartesian and cocartesian monoidal categories are (symmetric) purely monoidal categories.

In fact, suppose $\cat C$ is a category with finite products, and that for some objects $A$ and $B$  we have an isomorphism
\(
	A\times B\iso 1,
\)
where $1$ is the terminal object of $\cat C$ and $\times$ denotes the product. We show that $A$ is terminal. Since $A\times B\iso 1 $, it follows that $A\times B$ is itself terminal. Consequently, for every $X$ in $\cat C$, there exist maps $a\colon X\to A$ and $b\colon X\to B$ such that $(a,b)$ is the unique map from $X$ to $A\times B$. If $a'\colon X\to A$ is another map from $X$ to $A$, then $(a,b)=(a',b)$, as $A\times B$ is terminal. Hence, $a=a'$, and thus $A$ is terminal. 

We have therefore shown that all weakly invertible objects in a cartesian monoidal category are terminal, and hence uniquely isomorphic to the monoidal unit.
\end{example}
We establish now an important property of \emph{symmetric} purely monoidal categories; however, the proof relies on some results that will be introduced later in this paper.
\begin{proposition}
\label{const-pi}
$\catPu$ is a reflective subcategory of $\SMC$ (regarding both as 1-categories).
\end{proposition}
\begin{proof}
    In this proof we will explicitly exhibit a way to associate to any symmetric monoidal category $\cat M$ a symmetric purely monoidal category $\PF M$ and a symmetric monoidal functor $\prj M\colon \cat M\to\PF M$, since we will soon need this construction. However, we will not prove here that this yields in fact a reflection, as this will be a simple consequence of a later result (see \cref{theHTT} and \cref{RefCoref}). 
    
    The construction  we present here is based on a work by Vitale (\cite{VITALE}), where the existence of a certain notion of cokernel for morphisms between 2-groups is proven. We report here (an adaptation of) the construction of the cokernel in \cite{VITALE} applied to the inclusion of the Picard 2-group into a symmetric monoidal category, and we check that it gives rise to a purely monoidal category.  
    
    Let us first fix a symmetric monoidal category $\cat M$ and then proceed by steps. Recall that $\GP M$ is the Picard 2-group of $\cat M$ (see \cref{Pic-def}).
{
\setlength{\leftmarginii}{.7cm}
\begin{enumerate}
\item 
The first step is to define a bicategory $\bic M$ as follows.
\begin{itemize}
	\item The objects of $\bic M$ are the same as those of $\cat M$.
	\item Given two objects $X$ and $Y$ in $\cat M$, a 1-arrow in $\bic M$ from $X$ to $Y$ is given by a pair $(A, f)$, where $A$ is an object of $\PIC M$, and 
\(
	f\colon X\to Y\ten A
\)
is an arrow in $\cat M$.
	\item Given two 1-arrows in $\bic M$
\begin{equation*}
\begin{tikzcd}
	X\arrow[r,"{(A,f)}"]&Y\arrow[r, "{(B,g)}"]&Z,
\end{tikzcd}
\end{equation*}
the composition is given by $(B\ten A, h)$, where $h$ is given by the following composition of arrows in $\cat M$.
\begin{equation*}
\begin{tikzcd}[column sep = 3em]
	X\arrow[r, "f"] &Y\ten A\arrow[r, "g\ten A"]& (Z\ten B)\ten A\arrow[r, "{a_{Z,B,A}}"] &Z\ten(B\ten A).
\end{tikzcd}
\end{equation*}
\item The 1-identity on an object $X$ in $\bic M$ is given by $(I, r_X^{-1})$ ($r$ being the right unitor).
\item Given parallel 1-arrows $(A,f)$, $(A',f')\colon X\to Y$, we define a 2-arrow $\alpha\colon (A,f)\tto (A',f')$ to be just an arrow $\alpha\colon A\to A'$ in $\PIC M$ such that the following triangle commutes in $\cat M$.
\begin{equation}
\label{class-eq}
\begin{tikzcd}[column sep=1.5em, row sep=2.3em]
	&X\arrow[dl, "f"']\arrow[dr,"{f'}"]
\\
	Y\ten A\arrow[rr, "Y\ten\alpha"'] && Y\ten A'
\end{tikzcd}
\end{equation}
\item Vertical composition and 2-identities are just composition and identities in $\PIC M$, and horizontal composition is just the tensor product in $\PIC M$.
\item The unitors and associators of $\bic M$ (as a bicategory) are given by the unitors and associators of $\cat M$ (as a monoidal category).
\end{itemize}
This makes $\bic M$ into a bicategory where every 2-arrow is invertible (being an invertible arrow in $\PIC M$).
\item 
Having defined $\bic M$, we now define $\PF M$ as the classifying category of $\bic M$. 
Explicitly, we have that $\PF M$ is a category having
\begin{itemize}
	\item the same objects as $\bic M$ (and therefore the same objects as $\cat M$);
	\item as arrows, 2-isomorphism classes of 1-arrows of $\bic M$.
\end{itemize} 
We denote the 2-isomorphism class of a 1-arrow $(A,f)$ in $\bic M$ by $[A,f]$. 
By definition, given objects $X$ and $Y$ in $\cat M$, and 1-arrows 
\((A,f),\, (A',f')\colon X\to Y\)
in $\bic M$, we have that these determine equal arrows
\( [A,f]=[A',f']\colon X\to Y\)
in $\PF M$ if and only if there exists an isomorphism $\alpha\colon A\to A'$ in $\cat M$ such that the above triangle \ref{class-eq} commutes. 
\item
We have a functor
\[
	\prj M\colon \cat M\to \PF M
\]
mapping each object of $\cat M$ to itself in $\PF M$, and each morphism $f\colon X\to Y$ in $\cat M$ to $[I, \mcomp{f,r_Y^{-1}}]$:
\begin{equation*}
\begin{tikzcd}
	X\arrow[r, "f"]& Y\arrow[r, "r_Y^{-1}"] & Y\ten I.
\end{tikzcd}
\end{equation*}
\item
We now introduce a symmetric monoidal structure on $\PF M$ and observe that $\prj M$ is a symmetric monoidal functor.
\begin{itemize}
	\item The tensor product of objects is the same as in $\cat M$.
\item Given two arrows $[A,f]\colon X\to Y$ and $[A',f']\colon X'\to Y'$ in $\PF M$, their tensor product is given by $[A\ten A', g]$, where $g$ is the following composition of maps in $\cat M$.
\begin{equation*}
\begin{tikzcd}[row sep=2em, column sep= 6.5em]
	X\ten X'\arrow[r, "{f\ten f'}"] &(Y\ten A)\ten(Y'\ten A')\arrow[dl, phantom,""{coordinate, name = C}]\arrow[dl,
	rounded corners,
	to path={-- ([shift={(2ex,0)}]\tikztostart.east)
	|-([shift={(0,2.5ex)}]\tikztotarget.north)
	--(\tikztotarget)}]
	\\
	Y\ten\bigl((A\ten Y')\ten A'\bigr)\arrow[r, "{Y\ten(b_{A,Y'}\ten A')}" ] & Y\ten\bigl((Y'\ten A)\ten A'\bigr)\arrow[dl,
	rounded corners,
	to path={-- ([shift={(2ex,0)}]\tikztostart.east)
	|-([shift={(0,2.75ex)}]\tikztotarget.north)
	--(\tikztotarget)}]
	\\
	(Y\ten Y')\ten(A\ten A')
\end{tikzcd}
\end{equation*}
(We recall that $b$ is the symmetry, and the unnamed winding arrows are compositions of the appropriate components of the associator of $\cat M$.)
\item The associator, the unitors and the braiding of $\PF M$ are just the image under $\prj M$ of the corresponding maps of $\cat M$.
\end{itemize}
With this additional structure, we have made $\PF M$ into a symmetric monoidal category. The symmetry of $\cat M$ (and not just the braiding) is used to prove the naturality of the braiding of $\PF M$.
\item We now need to check that $\PF M$ is pure. We use a property that we state and prove separately in \cref{pure-iso} below, as it will be needed again later. 

In order to verify that $\PF M$ is pure we prove condition \ref{Pure-id} of \cref{Pure-eqdef}. Consider a weakly invertible object $A$ of $\cat M$. Then $[A,\ell_A^{-1}]$ is a morphism from $A$ to $I$ in $\PF M$, and it is also an isomorphism by \cref{pure-iso} (we recall that $\ell$ is the left unitor).

Now, suppose
\begin{equation*}
	[A,f],\,[B,g]\colon I\to I
\end{equation*}
are isomorphisms in $\PF M$. Then, again by \cref{pure-iso}, $f$ and $g$ are isomorphisms in $\cat M$. Call $\phi=\mcomp{{\ell_A}^{-1},f^{-1},g,\ell_B}$. By naturality of $\ell$ we have $\comp{f}{(I\ten\phi)}=g$, yielding $[A,f]=[B,g]$ (by definition).\qedhere
\end{enumerate}}
\end{proof}
\begin{lemma}
\label{pure-iso}
    Consider a symmetric monoidal category  $\cat M$ and the associated functor $\prj M\colon \cat M\to \PF M$ described in \cref{const-pi}.
    Then a morphism
\(
	[A, f]\colon X\to Y
\)
in $\PF M$ is an isomorphism if and only if $f$ is an isomorphism in $\cat M$.
In particular, $\prj M$ reflects isomorphisms.
\end{lemma}
\begin{proof}
As reported in \cite{VITALE}, if $f$ is an isomorphism in $\cat M$, and we consider an adjoint equivalence $(A,B,h,k)$ (see \cref{basic-2g}), then an inverse of $[A,f]$ in $\PF M$ is $[B, g]$, where $g$ is the composition of the following maps.
    \begin{equation*}
\begin{tikzcd}[row sep=2em]
	Y\arrow[r, "r_Y^{-1}"]& Y\ten I\arrow[r, "Y\ten h"] &Y\ten (A\ten B)\arrow[dll,
	rounded corners,
	to path={-- ([shift={(2ex,0)}]\tikztostart.east)
	|-([shift={(0,2.3ex)}]\tikztotarget.north)
	--(\tikztotarget)}]\arrow[d, phantom, "\scriptstyle{a_{Y,A,B}^{-1}}"{xshift=5.4em}]
\\
(Y\ten A)\ten B\arrow[rr, "f^{-1}\ten B"'] & &X\ten B
\end{tikzcd}
\end{equation*}

Vice versa, suppose $[A,f]$ is an isomorphism in $\PF M$, with inverse $[B, g]$.
Then, from the fact that 
\(
	\mcomp{{[A,f]},{[B,g]}}=\id X
\)
in $\PF M$ and by expanding the definition of composition and identity in $\PF M$, we obtain that in $\cat M$ the map $f$ has a left inverse and the map $g\ten A$ has a right inverse. Since $-\ten A$ is an equivalence (see \cref{basic-2g}), and equivalences preserve and reflect split monomorphisms and split epimorphisms, we conclude that $g$ has a right inverse as well.
Conversely, applying the same reasoning to the reverse composition $\mcomp{{[B,g]}, {[A,f]}}=\id Y$, we conclude that $g$ has a left inverse and $f$ has a right inverse. 
Since $f$ has both left and right inverses, it must be an isomorphism.

Finally, to see that $\prj M$ reflects isomorphisms consider an arrow $f\colon X\to Y $ in $\cat M$, and suppose that $\prj M(f)=[I,\mcomp{f, r_Y^{-1}}]$ is an isomorphism in $\PF M$.
Then, by the first part of this corollary, $\mcomp{f, r_Y^{-1}}$ is an isomorphism in $\cat M$, and therefore $f$ also is.
\end{proof} 
\begin{remark}
    As in \cite{VITALE}, the only step in the construction of \cref{const-pi} where symmetry, rather than just braiding, was used is in proving the naturality of the braiding of $\PF M$. It remains unclear whether this construction can be extended to braided monoidal categories.
\end{remark}
We conclude the section with another useful property related to the construction in \cref{const-pi}. This property corresponds exactly to Lemma~2.2 in \cite{VITALE}, except that the lemma in \cite{VITALE} is stated in the context of 2-groups; however, it  is easily shown to work in our context as well.
\begin{lemma}
\label{lemma-enrico}
Let $\cat M$ be a symmetric monoidal category, and $\prj M\colon \cat M\to \PF M$ the associated functor described in \cref{const-pi}.
Then, for any arrow $[A,f]\colon X\to Y$ in $\PF M$, the following diagram is commutative.
\begin{equation*}
\begin{tikzcd}[column sep= 6em]
	X\arrow[r, "{[A, f]}"]\arrow[d, "\prj M (f)"'] & Y
\\
Y\tens MA\arrow[r, "Y\ten_{\PF M}\mnhs A"'] &Y\tens M \I M\arrow[u, "\prj M(r_Y)"']
\end{tikzcd} 
\end{equation*}
Here, $\mnhs A= [A, \ell_A^{-1}]\colon A\to I$, with $\ell$ being the left unitor.
\end{lemma}
\section{Homotopy torsion theories}
\label{CNH}
{
\renewcommand{\nhs}{\Theta}
\renewcommand{\nh}[1]{\nhs({#1})}
\renewcommand{\ncomp}[3]{{#3}\bullet{#2}\bullet{#1}}
\renewcommand{\nncomp}[2]{{#2}\bullet{#1}}
\renewcommand{\theta}{\vartheta}
As mentioned in \intro{}, our intention is to develop a new torsion theory within the framework of homotopy torsion theories. In this section, we recall the essential notions needed for this purpose (see \cite{HTT} for more details).

We begin with the basic definition of \emph{category with nullhomotopies}, first introduced (in a slightly different form) by Grandis in \cite{GRANDIS97}.
\begin{definition}
\label{def-nh-struct}
Let $\cat C$ be any category. A \emph{nullhomotopy structure $\nhs$ on $\cat C$} is given by the following data.
\begin{enumerate}[1)]
	\item For every arrow $f$ in $\cat C$, a set $\Theta(f)$, called the \emph{set of nullhomotopies on $f$};
	\item for any triple of composable arrows 
\(
\begin{tikzcd}[cramped, sep=1.5em]
\cdot\arrow[r, "p"]&\cdot\arrow[r,"f"]&\cdot\arrow[r,"q"]&\cdot,
\end{tikzcd}
\)
a function 
\begin{equation*}
	\ncomp p - q\colon\nh f\to\nh{\mcomp{p,f,q}}
\end{equation*}
such that the following axioms are satisfied.
\begin{enumerate}[a)]
	\item \label[axiom]{nhidentity} Given $f\colon X\to Y$ and $\theta\in\nh f$, we have that 
\begin{equation*}
	\ncomp{\id X}\theta{\id Y}=\theta;
\end{equation*}
	\item \label[axiom]{nhassociativity}given a quintuple of composable maps
\(
\begin{tikzcd}[cramped, sep=1.5em]
	\cdot\arrow[r,"p'"]&\cdot\arrow[r, "p"]&\cdot\arrow[r, "f"]&\cdot\arrow[r, "q"]&\cdot\arrow[r,"q'"]&\cdot,
\end{tikzcd}
\) and given $\theta\in\nh f$, we then have that 
\[
\ncomp{(\comp {p'} {p})}{\theta}{(\comp q{q'})}=\ncomp{p'}{(\ncomp p{\theta}q)}{q'}.
\]
\end{enumerate}
\end{enumerate}

A \emph{category with nullhomotopies} is given by a pair $(\cat C, \nhs)$, where $\cat C$ is a category and $\nhs$ is a nullhomotopy structure on $\cat C$.

Given a category with nullhomotopies $(\cat C, \nhs)$, we say that a morphism $f$ in $\cat C$ is \emph{$\nhs$-trivial} (or \emph{nullhomotopic}) if $\nh f$ is non-empty. We say instead that an object in $\cat C$ is $\nhs$-trivial if its corresponding identity morphism is $\nhs$-trivial. Finally, given two objects $X$ and $Y$ in $\cat C$, we say that $X$ is \emph{$\nhs$-orthogonal to $Y$} if for all arrows $f\colon X\to Y$ we have that $\nh f$ is a singleton.
\end{definition}
In this section, within a category with nullhomotopies $(\cat C,\nhs)$, we will use the following notation to graphically represent a nullhomotopy $\theta$ on a map $f$,
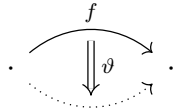
\begin{equation*}
\begin{tikzcd}[column sep =4em]
	\cdot & \cdot
	\arrow[""{name=0, anchor=center, inner sep=0}, "f", curve={height=-1.5em}, from=1-1, to=1-2]
	\arrow[""{name=1, anchor=center, inner sep=0}, curve={height=1.5em}, dotted, from=1-1, to=1-2]
	\arrow["\theta"{yshift=.1em, xshift=.1em}, shorten <=2pt, shorten=.3em, Rightarrow, from=0, to=1]
\end{tikzcd}
\end{equation*}
as if $\theta$ was a 2-cell in a 2-category from $f$ to a phantom 0-arrow (of course, in general, we may not be in a 2-category and we may not have a 0-arrow). Moreover, in light of \cref{nhidentity,nhassociativity} of \cref{def-nh-struct}, given maps $f\colon X\to Y$ and $g\colon Y\to Z$, and nullhomotopies $\phi\in\nh f$ and $\psi\in\nh g$, we write $\nncomp\phi g$ for $\ncomp{\id X}\phi g$ and $\nncomp f\psi$ for $\ncomp f\psi{\id Z}$.

Next, we introduce the natural notions of kernel and cokernel in this context, i.e.\ the notions of  \emph{homotopy kernel} and \emph{homotopy cokernel}.
\begin{definition}
\label{def-hker}
	Let $(\cat C,\nhs)$ be a category with nullhomotopies, and let $f\colon X\to Y$ be an arrow in $\cat C$. 
A \emph{homotopy kernel of $f$ with respect to $\nhs$} (or a \emph{$\nhs$-kernel of $f$}) is a triple $(K, k,\theta)$, where $k\colon K\to X$ is a morphism in $\cat C$ and $\theta\in\nh {\comp kf}$, with the property that for any other triple $(H,h,\phi)$, where $h\colon H\to X$ is a morphism in $\cat C$ and $\phi\in\nh{\comp hf}$, there exists a unique morphism $h'\colon H\to K$ such that $\comp {h'} k=h$ and $\nncomp {h'}{\theta}=\phi$.
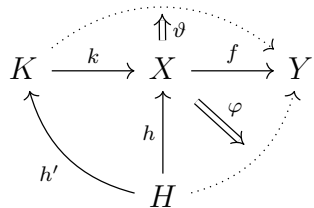
\begin{equation*}
\begin{tikzcd}
	K & X & Y \\
	& H
	\arrow["k", from=1-1, to=1-2]
	\arrow[""{name=0, anchor=center, inner sep=0}, curve={height=-2.5em}, dotted, from=1-1, to=1-3]
	\arrow["f", from=1-2, to=1-3]
	\arrow["{h'}", bend left, from=2-2, to=1-1]
	\arrow["h", from=2-2, to=1-2]
	\arrow[""{name=1, anchor=center, inner sep=0}, bend right, dotted, from=2-2, to=1-3]
	\arrow["\theta"{swap,yshift=-.1em,xshift=.1em},Rightarrow,from=1-2, to=0, shorten=.15 em]
	\arrow["\phi", shorten >=4pt, Rightarrow, from=1-2, to=1, shorten=.5 em]	
\end{tikzcd}
\end{equation*}
We will sometimes just write $K$ or $k$ to denote the $\nhs$-kernel $(K,k,\theta)$. 
Of course \emph{homotopy cokernels} are defined dually.

A \emph{homotopy exact sequence} in $(\cat C,\nhs)$ (or a \emph{$\nhs$-exact sequence}) will be given by a pair of composable morphisms
\(
\begin{tikzcd}[cramped, sep=1.5em]
X\arrow[r, "f"]&Y\arrow[r,"g"]&Z
\end{tikzcd}
\)
and a nullhomotopy $\theta\in\nh{\comp fg}$ such that $(X, f,\theta)$ is the $\nhs$-kernel of $g$ and $(Z,g,\theta)$ is the $\nhs$-cokernel of $f$.
\end{definition}

We are now ready to define homotopy torsion theories.
\begin{definition}
\label{HTTdef}
    Let $(\cat C,\nhs)$ be a category with nullhomotopies. A \emph{homotopy torsion theory} on $(\cat C,\nhs)$ (or a \emph{$\nhs$-torsion theory} on $\cat C$) is given by a pair $(\cat T,\cat F)$ of full and replete subcategories of $\cat C$ such that
    \begin{enumerate}
        \item \label[axiom]{HTTdef1} given objects $T$ in $\cat T$ and $F$ in $\cat F$, then $T$ is $\nhs$-orthogonal to $F$;
        \item \label[axiom]{HTTdef2} for every object $X$ in $\cat C$ there exists a $\nhs$-exact sequence
        \begin{equation*}
        \begin{tikzcd}
        {T_X} & X & {F_X},
	   \arrow["{t_X}"', from=1-1, to=1-2]
	   \arrow[""{name=0, anchor=center, inner sep=0}, bend left, dotted, from=1-1, to=1-3]
	   \arrow["{f_X}"', from=1-2, to=1-3]
	   \arrow["{\theta_X}"{yshift=-.2em}, shorten >=1pt, Rightarrow, from=1-2, to=0]
        \end{tikzcd}\end{equation*}
    with $T_X$ in $\cat T$ and $F_X$ in $\cat F$.
    \end{enumerate}
\end{definition}
}
\section{A nullhomotopy structure \texorpdfstring{$\nhs$}{Θ} for monoidal categories}
\label{Null-struct}
Let us now return to monoidal categories. In \cref{Pre,PMC} we considered two special classes of (symmetric) monoidal categories---namely 2-groups and purely monoidal categories---which could serve as the torsion and torsion-free objects in the new torsion theory we are developing. 
However, to construct a torsion theory, we first need a notion of trivial map and of exact sequence, and, naturally, we make use of the machinery just described in \cref{CNH}: our goal is now to endow the category of (symmetric) monoidal categories with a nullhomotopy structure, 
that is,  to define a set of nullhomotopies on any given (symmetric) monoidal functor. 

We know that between any two monoidal categories there always exists the constant functor at the unit object of the codomain. 
While this functor is not a zero (1-)morphism---since $\MC$ does not have a zero object---we will use it as a ``reference trivial morphism''. Specifically, we will consider any functor naturally isomorphic to it as trivial, and use these natural isomorphisms as nullhomotopies.
More precisely, we introduce the following definition. 
\begin{definition}
\label{MonNullStruct}
Given any monoidal functor 
\(
	F\colon\cat M\to\cat N,
\)
define the set
\begin{equation*}
	\nh F=\{\phi\colon F\tto\cst{\I N}\mid\phi\textnormal{ is a monoidal natural isomorphism}\}
\end{equation*}
(where $\cst{\I N}$ is the constant functor at $\I N$).
Given monoidal functors
\begin{equation*}
\begin{tikzcd}
	\cat M'\arrow[r,"P"]&\cat M\arrow[r,"F"] &\cat N\arrow[r,"Q"]&\cat N'
\end{tikzcd}
\end{equation*}
and given $\phi\in\nh F$, we define the composition
\begin{equation*}
\ncomp P\phi Q\in\nh{QFP}, 
\end{equation*}
having components
\begin{equation*}
	(\ncomp P\phi Q)_X=\comp{Q(\phi_{PX})}{e_Q^{-1}}
\end{equation*}
for all $X$ in $\cat M'$, where we recall that $e_Q$ is the coherence map of $Q$ defined in \cref{MonCat} (in other words, $\ncomp P\phi Q$ is the natural transformation obtained by vertically composing the horizontal composition $\mcomp{\id P,\phi,\id Q}$ with the constant natural transformation $e_Q^{-1}\colon\cst{(Q\I N)}\tto\cst{\I{ N'}}\colon \cat M'\to\cat N'$).
\end{definition}
\begin{remark}
    It is easy to verify that the definition of $\nhs$ we have just given makes $(\MC,\nhs)$ into a category with nullhomotopies. Clearly, if we restrict the definition of $\nhs$ to symmetric monoidal functors, we get a nullhomotopy structure on $\SMC$ (which we still call $\nhs$), making $(\SMC,\nhs)$ a category with nullhomotopies as well. 
\end{remark}
From now on, unless otherwise specified, $\nhs$ will always denote the specific nullhomotopy structure of \cref{MonNullStruct}.
\section{A \texorpdfstring{$\nhs$}{Θ}-torsion theory for symmetric monoidal categories}
Thus far, we have established the key ingredients for our 2-dimensional generalisation of the torsion theory on commutative monoids recalled in \intro{}, as summarised in the following table.
\begin{center}
\begin{tblr}{Q[4.3cm, valign=m, halign=r]Q[1.4cm, valign=m, halign=c]Q[4.3cm, valign=m, halign=l]}
\hline
  Commutative monoids & \(\implies\) & Symmetric monoidal categories \\\hline
  Trivial (=constant) maps    & \(\implies\) & $\nhs$-trivial (=naturally isomorphic to constant) functors\\\hline
  Abelian groups       & \(\implies\) & Symmetric 2-groups \\\hline
  Commutative pure monoids &\(\implies\) & Symmetric purely monoidal categories
  \\\hline
\end{tblr}   
\end{center}
We now need to verify that these ingredients do indeed yield the desired result. We do so in \cref{theHTT}, after proving the following preparatory lemma.
\begin{lemma}
\label{UniqueNH}
    Let $F\colon \cat M\to\cat N$ be a (symmetric) monoidal functor between (symmetric) monoidal categories. If $\cat N$ is pure, then $\nh F$ has at most one element. If in addition $\cat M$ is a 2-group, then $\nh F$ has exactly one element.
\end{lemma}
\begin{proof}
Suppose $\cat N$ is pure. If there exists a nullhomotopy $\phi\in\nh F$, then $\phi$ is unique. In fact, for each object $X$ in $\cat M$, the map $\phi_X$ is an isomorphism between $FX$ and $\I N$. Since $\I N$ is weakly invertible, so must $FX$ be. Uniqueness then follows from the fact that in a purely monoidal category, there exists exactly one isomorphism between any two weakly invertible objects.

If, in addition, we know that $\cat M$ is a 2-group, then for every object $X$ in $\cat M$, we have that $FX$ is weakly invertible in $\cat N$ (as $X$ is weakly invertible in $\cat M$); consequently, there is a unique isomorphism $\phi_X\colon FX\to\I N$. The maps $\phi_X$, for $X$ in $\cat M$, give rise to a monoidal natural isomorphism between $F$ and $\cst{\I N}$, where naturality and monoidality follow from the uniqueness of isomorphisms between weakly invertible objects in $\cat N$.
\end{proof}
\begin{proposition}
\label{theHTT}
    $(\catGr,\catPu)$ is a $\nhs$-torsion theory on $\SMC$.
\end{proposition}
\begin{proof}
We need to verify that \cref{HTTdef1,HTTdef2} of \cref{HTTdef} are satisfied.

\Cref{HTTdef1} ($\nhs$-orthogonality of symmetric 2-groups and symmetric purely monoidal categories) follows from \cref{UniqueNH}.
To prove \cref{HTTdef2}, fix a symmetric monoidal category $\cat M$ and consider the  sequence
\begin{equation}
\label{exactpresofM}
\begin{tikzcd}
	\GP M\arrow[r, "\inc M"] & \cat M \arrow[r, "\prj M"] &\PF M,
\end{tikzcd}
\end{equation}
where $\inc M \colon\GP M\to \cat M$ is the inclusion of the Picard 2-group of $\cat M$ in $\cat M$ itself (see \cref{Pic-def}), and $\prj M\colon \cat M\to\PF M$ is the functor constructed in the proof of \cref{const-pi}. 
We already know that $\GP M$ is a symmetric 2-group and that $\PF M$ is a symmetric purely monoidal category. 
Thus, it remains to check that \ref{exactpresofM} is a $\nhs$-exact sequence (see \cref{def-hker}). Call $\mnh M$ the unique nullhomotopy on $\mcomp{\inc M,\prj M}$. By uniqueness, its components are given, for any object $A$ in $\GP M$, by
\begin{equation*}
	(\mnh M)_A\equiv\mnhs A= [A, \ell_A^{-1}]\colon A\to I.
\end{equation*}
We now outline the proofs of the required universal properties for the homotopy exactness of \ref{exactpresofM}.

First, we check that $\bigl(\GP M, \inc M, \mnh M\bigr)$ is the $\nhs$-kernel of $\prj M$. Suppose we are given the diagram

\begin{equation*}
\begin{tikzcd}
{\GP M} & {\cat M} & {\PF M} \\
	& {\cat L}
	\arrow["{\inc M}", from=1-1, to=1-2]
	\arrow["{\prj M}", from=1-2, to=1-3]
	\arrow[""{name=0, anchor=center, inner sep=0}, curve={height=-40pt}, from=1-1, to=1-3]
	\arrow["F", from=2-2, to=1-2]
	\arrow[""{name=1, anchor=center, inner sep=0}, from=2-2, to=1-3, bend right]
	\arrow["{F'}", dashed, from=2-2, to=1-1, bend left]
	\arrow["{\mnh M}"'{xshift=.5ex, yshift=-.5ex}, shorten >=3pt, Rightarrow, from=1-2, to=0]
	\arrow["\phi", shorten =.7em, Rightarrow, from=1-2, to=1]	
\end{tikzcd}
\end{equation*}
where $F\colon\cat L\to\cat M$ is a symmetric monoidal functor, $\phi\in\nh {\mcomp{F,\prj M}}$ and the unnamed arrows are constant functors at the units of the respective codomains.
We need to prove that there exists a unique symmetric monoidal functor $F'\colon \cat L\to\GP M$ such that $\mcomp{F',\inc M}=F$ and $\nncomp {F'}{\mnh M}=\phi$. In fact, it suffices to show that $FX$ is weakly invertible for all objects $X$ in $\cat L$, and that $Ff$ is an isomorphism for all morphisms $f$ in $\cat L$, so that $F$ lifts to the Picard 2-group of its codomain.
Uniqueness then follows from the fact that $\inc M$ is an inclusion, and the condition on nullhomotopies is trivial because there can exist at most one nullhomotopy on any functor into a pure monoidal category by \cref{UniqueNH}.
With this in mind, let $X$ be an object in $\cat L$. Then $\phi_X$ is an isomorphism between $\prj M(F(X))=F(X)$ and $I$ in $\PF M$.
We thus have
\(
	\phi_X=[A, u]
\) for some weakly invertible object $A$ in $\cat M$ and some morphism $u\colon FX\to I\ten A$ in $\cat M$.
By \cref{pure-iso}, $u$ is actually an isomorphism in $\cat M$, yielding
\(
	FX\iso I\ten A\iso A.
\)
Given that $A$ is weakly invertible, $FX$ must be as well, as desired. Next, let  $f\colon X\to Y$ be any morphism in $\cat L$: by the naturality of $\phi$ we obtain that
\(
\comp{\prj M(F(f))}{\phi_Y}=\phi_X,
\)
 from which it follows that $\prj M(F(f))$ is an isomorphism. Since $\prj M$ reflects isomorphisms by \cref{pure-iso}, we conclude that $F(f)$ is indeed an isomorphism.

 It now remains to check that $\bigl(\PF M, \prj M, \mnh M\bigr)$ is the $\nhs$-cokernel of $\inc M$. Consider the following diagram.
 \begin{equation*}
\begin{tikzcd}
{\GP M} & {\cat M} & {\PF M} \\
	& {\cat N}
	\arrow["{\inc M}", from=1-1, to=1-2]
	\arrow["{\prj M}", from=1-2, to=1-3]
	\arrow[""{name=0, anchor=center, inner sep=0}, curve={height=-40pt}, from=1-1, to=1-3]
	\arrow["G", from=1-2, to=2-2]
	\arrow["{G'}", dashed, from=1-3, to=2-2, bend left]
	\arrow[""{name=1, anchor=center, inner sep=0}, from=1-1, to=2-2, bend right]
	\arrow["{\mnh M}"'{xshift=.5ex, yshift=-.5ex}, shorten >=3pt, Rightarrow, from=1-2, to=0]
	\arrow["\psi"', shorten=.7em, Rightarrow, from=1-2, to=1]	
\end{tikzcd}
\end{equation*}
Here, $G\colon \cat M \to\cat N$ is a symmetric monoidal functor with a nullhomotopy $\psi\in\nh{\mcomp{\inc M, G}}$; the unnamed arrows are the constant functor at the appropriate unit objects.
We need to show that there exists a unique symmetric monoidal functor $G'\colon \PF M\to \cat N$ such that $\mcomp{\prj M, G'}=G$ and $\nncomp {\mnh M}{G'}=\psi$. The first condition determines $G'$ uniquely on objects: given any $X$ in $\cat M$ we must have $G'X=G'(\prj M (X))=GX$. To determine the action of $G'$ on arrows, we rely on \cref{lemma-enrico}. 
Using this lemma together with the condition $\nncomp {\mnh M}{G'}=\psi$, it follows that for any arrow $[A,f]\colon X\to Y$ in $\PF M$, $G'([A,f])$ must  be the composition of the following arrows in $\cat N$.
\begin{equation*}
\begin{tikzcd}[column sep= 4.4em, row sep=2em]
	GX\arrow[r, "Gf"] 
	&G(Y\tens M A)\arrow[r, "(m_G)_{Y,A}^{-1}"]
	&GY\tens N GA \arrow[dll,
	rounded corners,
	to path={-- ([shift={(2ex,0)}]\tikztostart.east)
	|-([shift={(0,1em)}]\tikztotarget.north)
	--(\tikztotarget)}]\arrow[d, phantom,"\scriptstyle {GY\tens N \psi_A}"{near start, xshift=6em, yshift=-.5em}]
	\\
	GY\tens N \I N\arrow[rr, "(r_{\cat N})_{GY}"'].
	&&GY
\end{tikzcd}
\end{equation*}
(We also use the fact that the coherence maps $m_G$ and $m_{G'}$ must be equal, due to the identity $G=\mcomp{\prj M, G'}$.)
It is now routine to check that this uniquely defines $G'$ as a symmetric monoidal functor satisfying the required properties.
\end{proof}
We conclude the section with a few observations.
\begin{remark}
\label{RefCoref}
   By Corollary 7.8 in \cite{HTT} we have that, viewed as 1-categories, $\catGr
    $ is a coreflective subcategory of $\SMC$ and $\catPu$ is a reflective subcategory of $\SMC$, thus confirming \cref{2g-coref} and completing the proof of \cref{const-pi}.
\end{remark}
\begin{remark}
\label{2-dim-aspects}
    The homotopy torsion theory of \cref{theHTT} features some interesting 2-dimensional aspects: the homotopy kernels and homotopy cokernels involved in \cref{theHTT} can be shown to be particular instances of bilimits and bicolimits, respectively, and the adjunctions mentioned in \cref{RefCoref} can, in fact, be extended to biadjunctions. However, a fully 2-dimensional approach is still under investigation.
\end{remark}
\newcommand{\DMC}{\cats{Z}}
\begin{remark}
    The intersection $\catGr\cap\catPu$ of torsion and torsion-free objects in the homotopy torsion theory of \cref{theHTT} is given by the category $\DMC$ of (symmetric) monoidal categories that are monoidally equivalent to the terminal monoidal category $\catone$. Such categories are characterised by the fact that there exists exactly one morphism between any two of their objects. Despite the fact that $\DMC$ can be proven to be a reflective and coreflective subcategory of $\SMC$ (and of $\MC$), the $\nhs$-torsion theory $(\catGr,\catPu)$ is \emph{not} a pretorsion theory in the sense of \cite{FACCHINI}. A simple way to prove this is to observe that, unlike in the setting of pretorsion theories, $\nhs$-cokernels in this case are not epimorphisms,  as shown in the following counterexample.
\end{remark}
\begin{example}
Let $C=\{1,-1\}$ be the 2-element (multiplicative) group, and let $\cat M=\ioz(C)$ and $\cat N=\ioo(C)$ (with $\ioz$ and $\ioo$ defined in \cref{HomInv}).
Notice that there exists a unique (monoidal) functor $H$ from $\cat M$ to $\cat N$.
Since $\cat M$ is a 2-group, it easy to check that in the $\nhs$-presentation of $\cat M$
\begin{equation*}
\begin{tikzcd}[column sep=4em]
	\GP M\arrow[r, "\inc M"] &\cat M\arrow[r, "\prj M"]&\PF M,
\end{tikzcd}
\end{equation*}
we have that $\GP M=\cat M$, and the category $\PF M$ has two objects and exactly one arrow between any two of them.
There are exactly two functors from $\PF M$ to $\cat N$: one---call it $F$---which maps  the two non-identity isomorphisms of $\PF M$ to $1$ in $\cat N$, and the other---call it $G$---that sends them to $-1$ in $\cat N$. We have that $\comp {\prj M}F=\comp{\prj M}G=H$ and so ${\prj M}$ is not an epimorphism in $\SMC$.
\end{example}
\section{Homotopy invariants}
\label{HomInv}
{
In this section we discuss how the ``1-dimensional'' torsion theory on $\CMon$ described in \intro{} relates to the homotopy torsion theory we developed on $\SMC$. We consider three main ways to associate a commutative monoid to a symmetric monoidal category: by taking isomorphism classes, by taking connected components, and by considering the endomorphisms of the unit object. We now examine how each of these constructions interacts with the two aforementioned torsion theories.

We denote by $\MGPF$, $\MPFF$, $\minc{}$ and $\mprj{}$ the functors and natural transformations associated to the torsion theory in $\CMon$, forming for each commutative monoid $M$ the exact presentation
\[
\begin{tikzcd}
  \MGP M\arrow[r, "\minc M"] & M\arrow[r, "\mprj M"] &\MPF M,
\end{tikzcd}
\]
with $\MGP M$ an abelian group and $\MPF M$ a pure commutative monoid.

Consider first the functor
\[
	\K\colon\SMC\to\CMon
\]
that assigns to each (small) symmetric monoidal category \( \cat{M} \) the commutative monoid \( \K\cat{M} \), whose elements are isomorphism classes of objects in \( \cat{M} \), with the monoid operation induced by the tensor product of \( \cat{M} \) and the monoid unit given by the isomorphism class of the unit object of \( \cat{M} \).
One can readily verify that the functor $\K$ ``commutes'' with the functors $\GPF$ and $\PFF$,  yielding a ``homotopy-invariant'' exact sequence of commutative monoids associated to each symmetric monoidal category, as formalised in the following proposition.
\begin{proposition}
	Let $\cat M$ be a symmetric monoidal category. Then we have the following isomorphism of short exact sequences of monoids.
\begin{equation*}
\begin{tikzcd}
	\K\GPF\cat M\arrow[r, "\K\inc M"]\arrow[d, "\iso" description]&\K\cat M\arrow[r, "\K\prj M"]\arrow[d, equal]&\K\PFF \cat M\arrow[d, "\iso" description]
	\\
	\MGPF\K\cat M\arrow[r, "\minc{\K\cat M}"']&\K\cat M\arrow[r, "\mprj{\K\cat M}"']&\MPFF\K\cat M
\end{tikzcd}
\end{equation*}

If we call $\Pic M=\K\GPF\cat M\iso\MGPF\K\cat M$ (this is the standard Picard group of $\cat M$), and $\Pure M=\K\PFF \cat M\iso\MPFF\K\cat M$, then the short exact sequence of commutative monoids
\begin{equation*}
\begin{tikzcd}
	\Pic M\arrow[r] &\K\cat M\arrow[r] &\Pure M
\end{tikzcd}
\end{equation*}
is a ``homotopy invariant'' of $\cat M$, in the sense that symmetric monoidal equivalences induce isomorphisms between the corresponding exact sequences. 
\end{proposition}

Next, consider the following pairs of functors.
\begin{equation*}
\begin{tikzcd}[column sep =2em]
	\SMC\arrow[r,bend left, "\piz"]&\CMon\arrow[l, bend left, "\ioz"] &\SMC\arrow[r, bend left, "\pio"] &\CMon\arrow[l, "\ioo", bend left]
\end{tikzcd}
\end{equation*}
\begin{itemize}
	\item \(\piz\) assigns to each symmetric monoidal category \( \cat{M} \) the commutative monoid \( \piz(\cat{M}) \), whose elements are the connected components of \( \cat{M} \), with the monoid operation induced by the tensor product of \( \cat{M} \) and the monoid unit given by the connected component of the unit object of \( \cat{M} \).
\item
\(\ioz\) assigns to each commutative monoid $M$ the discrete monoidal category \( \ioz(M) \), whose objects are the elements of $M$, with the tensor product and unit induced from the monoid operation and unit of $M$.
\item 
\(\pio\) assigns to each symmetric monoidal category \( \cat{M} \) the commutative monoid \( \pio(\cat{M}) \), consisting of the endomorphisms of \( \I M \), with the monoid operation given by composition and the unit given by \( \id{\I M} \).
\item
\(\ioo\) assigns to each commutative monoid $M$ the symmetric monoidal category \( \ioo(M) \), which has a single object whose endomorphisms are given by the elements of $M$. Composition and tensor product in \( \ioo(M) \) are both simply given by the monoid operation in $M$.
\end{itemize}
Applying $\piz$ or $\pio$ to $\nhs$-exact sequence of the form \ref{exactpresofM} does not yield in general exact sequences in $\CMon$. However the following holds.
\begin{proposition}
{
\label{prop-hom-inv}
\newcommand{\emk}{\cat M_k}
For any commutative monoid $M$, denote by $\emk$ the symmetric monoidal category $\ioj k M$, for $k\in\{0,1\}$. We have the following natural isomorphisms of short exact sequences, for $k\in\{0,1\}$.
\begin{equation*}
\begin{tikzcd}[column sep = 8em,every label/.append style = {font = \normalsize}]
	\pij k\GPF\emk
		\arrow[d, "\iso" description]
		\arrow[r, "\pij k({\inc{}}_{\emk})"] 
	&\pij k\emk
		\arrow[d, "\iso" description]
		\arrow[r,"\pij k({\prj{}}_{\emk})"] 
	&\pij k \PFF\emk 
		\arrow[d, "\iso" description]
\\
	\MGPF M
		\arrow[r, "\minc M"']
	& M
		\arrow[r, "\mprj M"']
	&\MPFF M
\end{tikzcd}
\end{equation*} 
}%
Therefore, the short exact sequence associated to $M$ in the torsion theory on $\CMon$ can be recovered in two ways from homotopy exact sequences arising in the homotopy torsion theory on $\SMC$. Furthermore, for any commutative monoid $M$, the following equivalences hold: $M$ is a group if and only if $\ioz(M)$ is a 2-group, if and only if $\ioo(M)$ is a 2-group; similarly, $M$ is a pure monoid if and only if $\ioz(M)$ is a purely monoidal category, if and only if $\ioo(M)$ is a purely monoidal category.
\end{proposition}
}
\section{Another 2-dimensional torsion theory}
\label{another-tt}
\newcommand{\tenp}[2]{{#1}^{\ten{#2}}}%
\newcommand{\torobj}[1]{\textnormal{T}(\cat{#1})}%
\newcommand{\torincl}[1]{T_{\cat{#1}}}%
\newcommand{\torfobj}[1]{\Phi(\cat{#1})}%
\newcommand{\torfincl}[1]{F_{\cat{#1}}}%
In this work, we have developed a homotopy torsion theory on the 2-category of symmetric monoidal categories starting from a classical torsion theory in the 1-category of commutative monoids. In this final section, we aim to carry out a similar  process, but this time we start from the paradigmatic torsion theory on the category of abelian groups, given by torsion and torsion-free abelian groups.

The most natural approach to this generalisation is to work in the 2-category $\catGr$ of symmetric 2-groups, as a 2-dimensional counterpart to abelian groups, while maintaining the same nullhomotopy structure $\nhs$ of \cref{MonNullStruct}, though restricted to symmetric monoidal functors between symmetric 2-groups; we continue to refer to this restriction by $\nhs$.

Let us begin by introducing some preliminary definitions.
\begin{definition}
\label{tor-symgrp}
Given an object $X$ in a symmetric 2-group $\cat G$, we say that $X$ is a \emph{torsion object} if there exists a positive integer $n$ such that $\tenp X n\iso I$ (where the tensor power $\tenp X n$ is defined recursively by $\tenp X 1=X$ and $\tenp X {(m+1)}=X\ten \tenp X m$). 

We say that a symmetric 2-group $\cat G$ is \emph{torsion} if all of its objects are torsion; we say that a symmetric 2-group $\cat G$ is \emph{torsion-free} if for any torsion object $X$ in $\cat G$ there exists a unique isomorphism from $X$ to $I$ in $\cat G$.
\end{definition}
We claim in the following proposition that torsion and torsion-free symmetric 2-groups form a $\nhs$-torsion theory on the category on symmetric 2-groups.
\begin{proposition}
Torsion and torsion-free symmetric 2-groups (as defined in \cref{tor-symgrp})  form a homotopy torsion theory on $\catGr$ with respect to the nullhomotopy structure $\nhs$ defined at the beginning of the section.
\end{proposition}
\begin{proof}
We sketch the main steps of the proof.  
\begin{enumerate}[(1)]
    \item \label{unique-nh-t-tf} 
    Since the image of a torsion object under a symmetric monoidal functor  is clearly a torsion object again, it follows directly from the definitions of torsion and torsion-free symmetric 2-group that $\nh F$ is a singleton for every symmetric monoidal functor $F$ from a torsion symmetric 2-group to a torsion-free symmetric 2-group.
    \item\label{unique-nh-tf} Since an object isomorphic to the unit is necessarily a torsion object, it follows that in a torsion-free symmetric 2-group, there is at most one morphism from any object to the unit. Consequently, $\nh F$ is either empty or a singleton for any symmetric monoidal functor $F$ with arbitrary domain and a torsion-free symmetric 2-group as codomain.
    \item Let us fix a symmetric 2-group $\cat G$. We can consider the full subcategory $\torobj G$ of its torsion objects: clearly, $\torobj G$ is a torsion symmetric 2-group, and will serve as the torsion part of $\cat G$.
    \item Consider the inclusion $\torincl G\colon \torobj G\inclusion\cat G$. Since we work in the category of symmetric 2-groups, we can directly apply the cokernel construction from \cite{VITALE} to obtain the $\nhs$-cokernel $\torfincl G\colon \cat G\to \torfobj G$ of $\torincl G$ in $\catGr$. This construction is essentially the same as the one we used in \cref{const-pi} and is more extensively described in \cite{VITALE}. Here, we briefly recall the main features needed for our purposes.
    \begin{itemize}
        \item The category $\torfobj G$ has the same objects as $\cat G$.
        \item A morphism from $X$ to $Y$ in $\torfobj G$ is an equivalence class  $[A,f]$, where $A$ is an object in $\torobj G$ and $f\colon X\to Y\ten A$ is a morphism in $\cat G$; we have $[A,f]=[B,g]\colon X\to Y$ if and only if there exists a morphism $\psi\colon A\to B$ in $\cat G$ such that $g=\comp f{(Y\ten\psi)}$.
        \item The tensor product of objects in $\torfobj G$ is the same as that of $\cat G$.
    \end{itemize}
    \item We have thus obtained a sequence \begin{equation}\label{ex-pres-G}
        \begin{tikzcd}
            \torobj G\arrow[r,"\torincl G"]& \cat G \arrow[r, "\torfincl G"]&\torfobj G,
        \end{tikzcd}
    \end{equation} 
    where $\torfincl G$ is a $\nhs$-cokernel of $\torincl G$. To complete the proof, it remains to show that $\torfobj G$ is torsion-free and that $\torincl G$ is the $\nhs$-kernel of  $\torfincl G$. We rely on the following property (\cref{torsion-implies-torsion}).
    \item \label{torsion-implies-torsion} If an object $X$ is torsion in $\torfobj G$ then it is torsion in $\cat G$. In fact, if $X$ is torsion in $\torfobj G$, then we have a positive integer $n$ and a morphism $[A,f]\colon\tenp X n\to I$ in $\torfobj G$, for some map $f$ in $\cat G$ and some object $A$ in $\torobj G$. Therefore, by the definition of the arrows in $\torfobj G$, it follows that $\tenp X n\iso I\ten A\iso A$ in $\cat G$. Since $A$ is torsion in $\cat G$, there exists a positive integer $m$ such that $\tenp A m \iso I$ in $\cat G$. Consequently, $\tenp X {(nm)}\iso\tenp {(\tenp X n)}m\iso \tenp A m\iso I$, and so $X$ is torsion in $\cat G$.
    \item To check that $\torfobj G$ is torsion-free, suppose that $X$ is a torsion object in $\torfobj G$. By \cref{torsion-implies-torsion}, $X$ is torsion in $\cat G$, and so we can consider the isomorphism $[X,{\ell_X}^{-1}]\colon X\to I$ in $\torfobj G$. Given two isomorphisms $[A,f],[B,g]\colon X\to I$ in $\torfobj G$, we have $g=\mcomp{f,(I\ten\psi)}$, with $\psi=\mcomp{{\ell_A}^{-1},f^{-1},g,\ell_B}$, and therefore $[A,f]=[B,g]$.
    \item To verify that $\torincl G$ is the $\nhs$-kernel of $\torfincl G$, consider a 2-group $\cat H$ and a symmetric monoidal functor $F\colon \cat H\to \cat G$  such that $\mcomp{F,\torfincl G}$ is $\nhs$-trivial. We only need  to check that for all $X$ in $\cat H$ we have that $FX$ is torsion in $\cat G$, that is, that $F$ actually takes values in the full subcategory $\torobj G$ of $\cat G$ (in fact, the condition on nullhomotopies is automatically satisfied thanks to \cref{unique-nh-tf}, and uniqueness is ensured by the fact that $\torincl G$ is an inclusion). If  $\mcomp{F,\torfincl G}$ is $\nhs$-trivial, then for all $X$ in $\cat H$, we have that $FX\iso I$ in $\torfobj G$. By \cref{torsion-implies-torsion}, it follows that $FX$ is torsion in $\cat G$, completing the proof. \qedhere
\end{enumerate}
\end{proof}
The torsion theory just introduced shares the same 2-dimensional properties mentioned in \cref{2-dim-aspects}. It is also well-behaved with respect to the functors $\piz$ and $\pio$ (see \cref{HomInv}), meaning that by applying these functors to the $\nhs$-exact sequence \ref{ex-pres-G} we obtain exact sequences of abelian groups corresponding to two different torsion theories. In particular, applying $\piz$ we get the classical presentation of the abelian group $\piz(\cat G)$ as an extension of its (classical) torsion-free part by its torsion subgroup; applying $\pio$ we get instead the trivial exact sequence 
\(
\begin{tikzcd}[cramped, sep=2em]
\pio(\cat G)\arrow[r,equal]&\pio(\cat G)\arrow[r]&0.
\end{tikzcd}
\)
\section*{Acknowledgements}
\addcontentsline{toc}{section}{Acknowledgements}
This research was partially conducted in affiliation with INdAM – Istituto Nazionale di Alta Matematica “Francesco Severi”, Gruppo
Nazionale per le Strutture Algebriche, Geometriche e le loro Applicazioni. The author sincerely thanks the reviewer for their careful reading and valuable feedback. Additionally, the author is deeply grateful to Sandra Mantovani for many meaningful discussions and insightful perspectives, and to Enrico Vitale for his helpful suggestions and remarks.

\bibliographystyle{elsarticle-num} 

\end{document}